\documentclass[11pt]{article}
\usepackage{amsfonts,amssymb,amsmath}

\usepackage[english]{babel}
\newcommand{\halmos}{\rule{1ex}{1.4ex}}

\newtheorem{ittheorem}{Theorem}
\newtheorem{itlemma}{Lemma}
\newtheorem{itproposition}{Proposition}
\newtheorem{itdefinition}{Definition}
\newtheorem{itremark}{Remark}

\newenvironment{theorem}{\addtocounter{equation}{1}
\begin{ittheorem}}{\end{ittheorem}}

\newenvironment{lemma}{\addtocounter{equation}{1}
\begin{itlemma}}{\end{itlemma}}

\newenvironment{proposition}{\addtocounter{equation}{1}
\begin{itproposition}}{\end{itproposition}}

\newenvironment{definition}{\addtocounter{equation}{1}
\begin{itdefinition}}{\end{itdefinition}}

\newenvironment{remark}{\addtocounter{equation}{1}
\begin{itremark}}{\end{itremark}}

\newenvironment{corollary}{\addtocounter{equation}{1}
\begin{itcorollary}}{\end{itcorollary}}

\newenvironment{proof}{\noindent {\em Proof}.\,\,\,}
{\hspace*{\fill}$\halmos$\medskip}

\newcommand{\beq}{\begin{eqnarray}}
\newcommand{\eeq}{\end{eqnarray}}

\newcommand{\beqt}{\begin{eqnarray*}}
\newcommand{\eeqt}{\end{eqnarray*}}

\newcommand{\be}{\begin{equation}}
\newcommand{\ee}{\end{equation}}

\newcommand{\bl}{\begin{lemma}}
\newcommand{\el}{\end{lemma}}

\newcommand{\br}{\begin{remark}}
\newcommand{\er}{\end{remark}}

\newcommand{\bt}{\begin{theorem}}
\newcommand{\et}{\end{theorem}}

\newcommand{\bd}{\begin{definition}}
\newcommand{\ed}{\end{definition}}

\newcommand{\bp}{\begin{proposition}}
\newcommand{\ep}{\end{proposition}}

\newcommand{\bc}{\begin{corollary}}
\newcommand{\ec}{\end{corollary}}

\newcommand{\bpr}{\begin{proof}}
\newcommand{\epr}{\end{proof}}

\newcommand{\bi}{\begin{itemize}}
\newcommand{\ei}{\end{itemize}}

\newcommand{\ben}{\begin{enumerate}}
\newcommand{\een}{\end{enumerate}}

\newcommand{\Z}{\mathbb Z}
\newcommand{\R}{\mathbb R}
\newcommand{\N}{\mathbb N}
\newcommand{\M}{\mathbb M}

\newcommand{\E}{\mathbb E}

\newcommand{\pee}{\mathbb P}

\newcommand{\s}{\ensuremath{\mathcal{S}}}

\newcommand{\om}{\ensuremath{\omega}}
\newcommand{\Om}{\ensuremath{\Omega}}

\newcommand{\si}{\ensuremath{\sigma}}
\newcommand{\eps}{\ensuremath{\epsilon}}

\begin{document}

\title{{\bf  Random walks on FKG-horizontally oriented lattices}}

\author{Nadine Guillotin-Plantard\footnote{Universit\'e Claude Bernard Lyon I,
LaPCS,
50, av. Tony-Garnier,
Domaine de Gerland,
69366 Lyon Cedex 07, France. E-mail: nadine.guillotin@univ-lyon1.fr}, Arnaud Le Ny\footnote{Eurandom, L.G. 1.48, TU Eindhoven, Postbus 513, 5600 MB Eindhoven, The Netherlands. E-Mail: leny@eurandom.tue.nl}}

\maketitle {\bf Keywords}: Markov chains, random environments, oriented graphs,
associated random variables,
recurrence vs transience, functional limit theorems.


\footnotesize
\begin{center}
{\bf Abstract:}
\end{center}
We study the asymptotic behavior of the simple random walk on oriented versions of $\Z^2$. The considered lattices are not directed  on the vertical axis but unidirectional on the horizontal one, with symmetric random orientations which are positively correlated. We prove that the simple random walk is transient and also prove a functional limit theorem in the space $\mathcal{D}([0,\infty[,\R^2)$ of c\`adl\`ag functions, with an unconventional normalization.

\normalsize

\vspace{12pt}
\newpage
\section{Introduction}

The use  of random walks as a tool in
mathematical physics is now well established and they have been for
example widely used in classical statistical mechanics to study critical phenomena
 (see \cite{FFS}). It has been recently observed that analogous methods in
  quantum statistical mechanics require the study of random walks on oriented lattices,
  due to the intrinsic non commutative character of the (quantum) world (see e.g. \cite{CP2,LR}).
   Although random walks in random and non-random environments have been intensively studied for
    many years, only a few  results on random walks on oriented lattices are known. The
      recurrence versus transience properties of simple random walks on oriented
       versions of $\Z^2$ are studied in \cite{CP} when the horizontal lines are
       unidirectional towards a random or deterministic direction. The interesting
       behavior of this model is that, depending on the orientation, the walk could
       be either recurrent or transient. In the deterministic "alternate" case, for
        which the orientations of horizontal lines are  alternated, i.e. oriented rightwards at one level and leftwards at the following level, the recurrence of the simple random walk is proved, whereas the transience
         naturally arises when the orientations are all identical in infinite regions.
         More surprisingly, it is also proved that the recurrent character of the simple
          random walk on $\Z^2$ is lost when the orientations are i.i.d. with zero mean.

In this paper, we prove that the transience of the simple random
walk still holds when the orientations are  symmetric and
positively correlated with a summable power law decay of
correlations. We also prove a functional limit
theorem for this walk with an unconventional normalization due to
the random character of the environment of the walk, solving an
open question of \cite{CP}. Our paper is organized as follows: the
description of  our model and the   results are stated in Section
2. Section 3 is devoted to the proofs while  illustrative
examples of orientations, coming from statistical mechanics, are
given in Section 4.

\section{Model and results}
\subsection{FKG-horizontally oriented lattices}
We consider  a canonical probability space $(\Om,\mathcal{B}, \pee )$ on which all the random
variables are defined, and denote  $\mathbb{E}$ (resp. {\rm Cov}) the expectation (resp. covariance) under $\pee$.  By {\em orientations},
we mean  a stationary symmetric family of $\{-1,+1\}$-valued  random
variables $(\eps_y)_{y \in \Z}$, with the following properties:
\ben
\item{Associated random variables:}\\
 For any $m \geq 0$, for any finite collection $(\eps_0,\dots,\eps_m)$, for any coordinatewise
 nondecreasing  functions $f,g$ on $\{-1,+1\}^m$,
\[
{\rm Cov} \big[ f(\eps_0,\dots,\eps_m);
g(\eps_0,\dots,\eps_m)\big] \geq 0.
\]
\item{Summable power-law decay of correlations:}\\
There exists $\alpha > 1$ such that
\[
\E \big[\eps_0 \eps_y]= \mathcal{O} \big( {\mid y \mid}^{-\alpha} \big) \; \rm{when} \;
\mid y \mid \longrightarrow \; + \infty.
\]
\een
In our set-up, these orientations are natural extensions of Rademacher random variables
of \cite{CP}.
 They have the same one-dimensional law ($\pee [\eps_0=+1]=\pee
 [\eps_0=-1]=\frac{1}{2}$) but are {\em not necessarily} independent.
 The notion of associated random variables (see \cite{NW})  is very natural in the context of Gibbs measures in
 statistical mechanics where it is equivalent to the FKG property of the joint distribution
 $\nu$ of the random field $\eps=(\eps_y)_{y \in \Z}$
 (\cite{FKG}). In such cases, we also  say that the orientations are {\em positively correlated}. Examples of such distributions are  ferromagnetic, possibly long range,  Ising models, described in Section 4 at the end of this paper.

We use these associated random variables to build our {\em FKG-horizontally oriented} lattices. These lattices are oriented versions of $\Z^2$: the vertical lines are  not oriented and the horizontal ones are unidirectional, the orientation at a level $y \in \Z$ being given by the random variable $\eps_y$ (say right if the value is $+1$ and  left if
 it is $-1$). More formally we give the
\bd[FKG-horizontally oriented lattices]
Let $\eps=(\eps_y)_{y \in \Z}$ be a sequence of $\{-1,+1\}$-valued, associated and symmetric random variables. The {\em FKG-horizontally oriented lattice} $\mathbb{L}^\eps=(\mathbb{V},\mathbb{A}^\eps)$ is the directed graph with vertex set $\mathbb{V}=\Z^2$ and edge set $\mathbb{A}^\eps$ defined by the condition that for $u=(u_1,u_2), v=(v_1,v_2) \in \Z^2$, $(u,v) \in \mathbb{A}^\eps$ if and only if
\ben
\item either $v_1=u_1$ and $v_2=u_2 \pm 1$
\item or $v_2=u_2$ and $v_1=u_1+ \eps_{u_2}$.
\een
\ed

\subsection{Simple random walk on $\mathbb{L}^\eps$}

We consider the usual simple random walk  $M=(M_n)_{n \in \mathbb{N}}$ on $\mathbb{L}^\eps$. It is a Markov chain with transition probabilities defined for all $(u,v) \in \mathbb{A}^\eps$ by
\[\pee[M_{n+1}=v  | M_n=u]=\;\left\{
\begin{array}{lll} \frac{1}{3}  \; &\rm{if} \;  (u,v)   \in \mathbb{A}^\eps&\\
\\
0 \; \; &\rm{otherwise.}&
\end{array}
\right.
\]
Its transience is proved in \cite{CP} for almost every orientation when the directions are given by a sequence of i.i.d. random variables $(\eps_y)_{y \in \Z}$, i.e. when the law $\nu$ of the random field $\eps$ is a product probability measure. We generalize this result in this positively correlated and possibly non-independent  context.

\bt \label{thm1} For $\nu$-a.e. realization of the orientation
$\eps$, the simple random walk on the FKG-horizontally oriented
lattice $\mathbb{L}^\eps$ is transient. \et

We also answer in this general set-up to an open question of
\cite{CP} and obtain a functional limit theorem with a suitable
and unconventional normalization. We consider a Brownian motion
$(W_{t})_{t\ge 0}$ and denote $(L_{t}(x))_{t \ge 0}$ its
corresponding local time at $x\in\R$. Moreover, we introduce a
pair of independent Brownian motions $Z_{+}(x), Z_{-}(x), x\geq
0$. We assume these processes to be defined on one probability
space and to be independent of each other so that the following
process is well-defined for all $t \geq 0$:
\begin{equation}\label{th}
\Delta_{t}=\int_{0}^{\infty}L_{t}(x)dZ_{+}(x)+\int_{0}^{\infty}L_{t}(-x)dZ_{-}(x).
\end{equation}
This process is a particular example from a family of processes obtained in \cite{KS} as functional limits of
$\Z$-valued  random walks in random sceneries. Moreover, it
has a continuous version which is self-similar with index
$\frac{3}{4}$ and has stationary increments. We also introduce a
real constant $m=\frac{1}{2}$, defined later as the mean of some geometric
random variables related to the  behavior of the walk in the horizontal direction.
\bt\label{thm2} The following
convergence holds: \be \label{flt} \Big(\frac{1}{n^{3/4}} M_{[nt]}
\Big)_{t \geq 0} \; \stackrel{\mathcal{D}}{\Longrightarrow}
\frac{m}{(1+m)^{3/4}}( \Delta_t,0)_{t \geq 0} \ee where
$\stackrel{\mathcal{D}}{\Longrightarrow}$  stands for 
convergence  in the space of c\`adl\`ag functions
$\mathcal{D}([0,\infty),\R^2)$ endowed with the Skorohod topology.
\et

\section{Proofs}
\subsection{Vertical and horizontal embeddings of the simple random walk}
 The simple random walk $M$  defined  on $(\Om,\mathcal{B},\pee)$ can be decomposed into a vertical and an horizontal part by restriction to the corresponding  axis. The vertical part is a simple random walk $Y=(Y_n)_{n \in \N}$ on the line. The (independent) $\si$-algebras generated by this vertical walk $Y$ and the orientation $\eps$ are denoted respectively by
\[
\mathcal{F}=\si(Y_n, \; n \in \N) \; \; \; {\rm and} \; \; \; \mathcal{G}=\si(\eps_y, \; y \in \Z).
\]
We also define for all $n \in \N$ and $y \in \Z$ the {\em local time} at level y of the walk $Y$ to be
\[
\eta_n(y)=\sum_{k=0}^n \mathbf{1}_{Y_k=y}.
\]
The horizontal embedding is a random walk with  $\mathbb{N}$-valued geometric jumps. More formally, a doubly infinite family $(\xi_i^{(y)})_{i \in \mathbb{N}^*, y \in \Z}$ of independent geometric random variables of parameter $p=\frac{1}{3}$ (and mean $m=\frac{1}{2}$) is given and one defines the embedded horizontal random walk $X=(X_n)_{n\in\N}$ by $X_0=0$ and for $n \geq 1$,
\[
X_n=\sum_{y \in \Z} \eps_y \sum_{i=1}^{\eta_{n-1}(y)} \xi_i^{(y)}
\]
with the convention that the last sum is zero when $\eta_{n-1}(y)=0$.
Of course, the walk $\M_n$ does not coincide with $(X_n,Y_n)$ but these objects are closely related: define for all $n \in \N$
\[
T_n=n + \sum_{y \in \Z} \sum_{i=1}^{\eta_{n-1}(y)} \xi_i^{(y)}
\]
to be the instant just after the random walk $M$ has performed its n$^{\rm{th}}$ vertical move. The following Lemma is proved in \cite{CP}.
\bl \label{lem1}
\ben
\item $M_{T_n}=(X_n,Y_n)$.
\item For a given orientation $\eps$, the transience of $(M_{T_n})_{n \in \N}$ implies the transience of $(M_n)_{n \in \N}$.
\een
\el

\subsection{Associated random variables}
The extension from  the i.i.d. case to our case is made possible by  a comparison
of the joint characteristic function of associated random variables with the product
 of the marginal ones, due to Newman {\em et al.} (\cite{NW}).

\bl \label{lem2}
Let $\eps=(\eps_y)_{y \in \Z}$ be a sequence of associated random variables. Then, for all $t \in \R, \; n \in \N$,
\be \label{Ass}
\Big| \E \Big[ e^{it\sum_{y \in \Z} \eps_y \eta_n(y)} \big| \mathcal{F} \Big]- \prod_{y \in \Z} \E\Big[e^{it \eps_y \eta_n(y)} \big| \mathcal{F} \Big] \Big| \leq \frac{1}{2}  t^2 \sum_{x \neq y} \eta_n(x) \eta_n(y) \E\big[ \eps_x \eps_y\big].
\ee
\el

\bpr It is based on Theorem 1 in \cite{NW}, which  states that for
a finite family of $p$ associated r.v.'s $(Z_1, \dots, Z_p)$ and
real numbers $(r_1,\dots,r_p)$, \be \label{NM} \Big|
\E\big[e^{i\sum_{k=1}^p r_k Z_k}\big] - \prod_{k=1}^p
\E\big[e^{ir_k Z_k}\big] \Big|\leq \frac{1}{2} \sum_{1\leq j\neq k
\leq p} |r_j||r_k| \rm{Cov}(Z_j,Z_k). \ee The sum and product of
the l.h.s of (\ref{Ass}) have a finite number of terms because
$\eta_n(y)=0$ for $|y| > n$. It is thus straightforward to derive
(\ref{Ass}) from (\ref{NM}) using the $\mathcal{F}$-measurability
of the local times $\eta_n(y)$, the associativity of $\eps$ and
its independence with the vertical walk $Y$. \epr
\subsection{Proof of the transience of the simple random walk}

The vertical walk $Y$ is  known to be recurrent and its asymptotic behavior is rather well controlled. The transience is due to the behavior of the embedded horizontal random walk $X$ and to exploit it we introduce a partition of $\Omega$ between typical or untypical paths of $Y$.

In all this proof, for any $i \in \N$, $\delta_i$ is a strictly positive  real number and we write $d_{n,i}=n^{\frac{1}{2}+\delta_i}$. Define the sets
\[
A_n=\big\{ \om \in \Om; \max_{0 \leq k \leq 2n} \; |Y_k| < n^{\frac{1}{2} +\delta_1} \big\} \; \cap \; \big\{ \om \in \Om; \max_{y \in \Z} \; \eta_{2n-1}(y) < n^{\frac{1}{2} +\delta_2}\big\}
\]
and
\[
B_n=\big\{\om \in A_n; \Big| \sum_{y \in \Z} \eps_y \eta_{2n-1}(y) \Big| > n^{\frac{1}{2} +\delta_3}\big\}.
\]
By Lemma \ref{lem1}, the transience of $M$ will be insured as soon as

\be \label{eqn1}
\sum_{n \in \N} \pee[X_{2n}=0;Y_{2n}=0] \; < \; \infty
\ee
and to do so we first decompose $\pee[X_{2n}=0;Y_{2n}=0]$ into
\be \label{pn123}
\pee[X_{2n}=0;Y_{2n}=0;A_n^c] + \pee[X_{2n}=0;Y_{2n}=0;B_n]
 + \pee[X_{2n}=0;Y_{2n}=0; A_n \setminus B_n].
\ee
 Some results of the i.i.d. case of \cite{CP} still hold and in particular we can prove using standard techniques the following
\bl \label{lem3}
For any $\delta_{1},\delta_2 > 0$,
\[
\sum_{n \in \N} \pee[X_{2n}=0;Y_{2n}=0;A_n^c] \; < \; \infty.
\]
\el
The second term of (\ref{pn123}) is also a generic term of convergent series due to the untypical character of the paths in $B_n$. Again from \cite{CP} with standard techniques, we have the
\bl \label{lem4}
For any $\delta_3>0$,
\[
\sum_{n \in \N} \pee[X_{2n}=0;Y_{2n}=0;B_n] \; < \; \infty.
\]
\el
Now, we denote
\[p_n=\pee[X_{2n}=0;Y_{2n}=0;A_n \setminus B_n].
\]
To prove the theorem, it remains to show that for some
$\delta_1,\delta_2,\delta_3>0$ \be \label{eqn2} \sum_{n \in \N}
p_n \; < \; \infty. \ee Decompose \be \label{pn} p_n=\E\Big[
\mathbf{1}_{Y_{2n=0}} \E \big[\mathbf{1}_{X_{2n=0}} \E
\big[\mathbf{1}_{A_n \setminus B_n}\big| \mathcal{F} \vee
\mathcal{G} \big] \big| \mathcal{F} \big] \Big]. \ee It is well
known that for the simple random walk $Y$, there exists $C>0$ s.t.
\be \label{srw} \pee[Y_{2n}=0] \sim C \cdot n^{-\frac{1}{2}}, \; n
\rightarrow + \infty \ee and we can prove  as in \cite{CP} the \bl
\label{lem5}

On the set $A_n \setminus B_n$, we have,
\be \label{sqrtln}
\pee \big[X_{2n}=0 \big| \mathcal{F} \vee \mathcal{G}\big]=
\mathcal{O} \Big( \sqrt{\frac{\ln{n}}{n}}\Big).
\ee

\el Hence, the transience of the simple random walk is a direct
consequence of the following

\bp \label{prop1} For  $\alpha >1$, it is possible to choose  $\delta_1,\delta_2,\delta_3>0$ such that there exists $\delta>0$ and
 \be \label{eqn3} \pee\big[ A_n \setminus
B_n \big| \mathcal{F} \big]= \mathcal{O} \big(n^{-\delta}). \ee \ep

\bpr
We first use an auxiliary centered Gaussian random variable with variance $d_{n,3}^2$, the Anderson's inequality and Plancherel's formula to get
\be \label{eqn4}
\pee\big[ A_n \setminus B_n \big| \mathcal{F} \big] \leq C \cdot n^{\frac{1}{2} + \delta_3} \cdot I_n
\ee
where
\[
I_n=\int_{-\pi}^{\pi} \E \big[e^{it \sum_{y \in \Z} \eps_y \eta_{2n-1}(y)} \big| \mathcal{F} \big] e^{-t^2 d_{n,3}^2/2} dt.
\]
To use that for $td_{n,3}$ small enough, $e^{-t^2 d_{n,3}^2/2}$
dominates the term under the expectation, we split the integral in
two parts. For $b_n=\frac{n^{\delta_2}}{d_{n,3}}$, we write
\[
I_n=I_n^{1} + I_n^2
\]
with
\beqt
I_n^1=\int_{|t| \leq b_n} \E \big[e^{it \sum_{y \in \Z} \eps_y \eta_{2n-1}(y)} \big| \mathcal{F} \big] e^{-t^2 d_{n,3}^2/2} dt\\
I_n^2=\int_{|t| > b_n} \E \big[e^{it \sum_{y \in \Z} \eps_y \eta_{2n-1}(y)} \big| \mathcal{F} \big] e^{-t^2 d_{n,3}^2/2} dt.
\eeqt
To control the integral $I_n^2$, we write
\beqt
|I_n^2| &\leq& C \int_{|t| > b_n} e^{-t^2 d_{n,3}^2/2} dt= \frac{C}{d_{n,3}} 
\int_{|s| > n^{\delta_2}} e^{-s^2/2} ds\\
&\leq& 2 \/ \frac{C}{d_{n,3}} \/  n^{-\delta_2} \/  e^{-n^{2 \delta_2} /2}
\eeqt
to get 
\[
|I_n^2|=\mathcal{O} \big(e^{-n^{2 \delta_2} / 2}).
\]

Let $(\eps'_y)_{y \in \Z}$ be a sequence of i.i.d. random variables with marginal distribution $\pee[\eps'_y=-1]=\pee[\eps'_y=+1]=\frac{1}{2}$ and denote  $I_n^0$ the integral which corresponds to $I_n^1$ in this case. Factorization is possible by independence and we write
\[
I_n^0=\int_{|t| \leq b_n}  \prod_{y \in \Z} \E \big[e^{it\eps'_y
\eta_{2n-1}(y)} \big| \mathcal{F} \big]e^{-t^2 d_{n,3}^2/2}
dt=\int_{|t| \leq b_n} \prod_{y \in \Z} \cos(\eta_{2n-1}(y)t)
e^{-t^2 d_{n,3}^2/2} dt
\]
and decompose
\[
I_n^1=I_n^0 + (I_n^1 - I_n^0).
\]
In order to get a validity of our result for any summable power law decay of correlations, we  estimate $I_n^0$ by  the following
\bl \label{lem6} For $\delta_3 > 2 \delta_2 $,
\[
|I_n^0|=\mathcal{O} \big(n^{-\frac{3}{4}+\frac{\delta_1}{2}}\big).
\]
\el
\bpr
We first use H\"older's inequality to get
\beqt
|I_n^0| &\leq& \prod_y \Big[ \Big( \int_{|t| \leq b_n} | \cos(\eta_{2n-1}(y)t)|^{\frac{2n}{\eta_{2n-1}(y)}} dt \Big)^{\frac{\eta_{2n-1}(y)}{2n}} \Big].
\eeqt
Denote for all $y \in \Z, n \in \N$, $p_{n,y}=\frac{\eta_{2n-1}(y)}{2n}$,
$C_n=\{y:\eta_{2n-1}(y)\ne 0\}$ and, for $y \in C_n$
\[J_{n,y}=
\int_{|t|\leq b_n} |\cos(\eta_{2n-1}(y) t)|^{1/p_{n,y}} dt=
\frac{1}{\eta_{2n-1}(y)}
\int_{|v|\leq b_n \eta_{2n-1}(y)} |\cos(v)|^{1/p_{n,y}} dv,\]
to get
$|I_n^0| \leq \prod_y J_{n,y}^{p_{n,y}}$. Now, using the fact that we work on $A_n$, we  choose $\delta_3>2 \delta_2$ in order to have
 $b_n\eta_{2n-1}(y)\rightarrow 0$ uniformly in $y$ when $n$ goes to infinity. Using the equivalence around zero between the cosinus and the exponential, one has
\begin{eqnarray*}
|I_n^0|&\leq& \prod_y \Big(\sqrt{\frac{2\pi}{2n\eta_{2n-1}(y)}} \; \Big)^{p_{n,y}}\\
& = & (2\pi)^{\frac{1}{2} \sum_y p_{n,y}}
\exp \big(-\frac{1}{2}\sum_y p_{n,y}\log(2n\eta_{2n-1}(y))\big).
\end{eqnarray*}
The vector $\textbf{p}=(p_{n,y})_{y \in C_n}$ defines a probability measure on  $C_n$ and we have
\[-\frac{1}{2}\sum_y p_{n,y}\log(2n\eta_{2n-1}(y))=
-\log 2n -\frac{1}{2}\sum_y p_{n,y}\log p_{n,y}=
-\log 2n+ \frac{1}{2}H(\textbf{p})\]
 where $H(\cdot)$ is the entropy of the probability vector
$\textbf{p}$, always bounded by $\log\textrm{card}C_n$. We thus have on the set $A_n$,
\[|I_n^0|\leq \sqrt{2\pi}\exp(-\log 2n + \frac{1}{2}\log (2 d_{n,1}))=
\frac{\sqrt{\pi d_{n,1}}}{n}=\mathcal{O}(n^{-\frac{3}{4}+\frac{\delta_1}{2}}).\]
\epr

From Lemma \ref{lem2}, it is possible  to compare $I_n^1$ to $I_n^0$ and control their difference.
\bl \label{lem7} For $\delta_3 > 2 \delta_2$ and $\beta=3 \delta_3 + \alpha -1 - 4 \delta_2$,
\be \label{beta}
|I_n^1-I_n^0| = \mathcal{O} \big(n^{- \beta}).
\ee
\el
\bpr
We have
\[
\Big| I_n^1 - I_n^0 \big| \leq J_n:=\int_{|t| \leq b_n} \Big| \E \big[e^{it \sum_{y \in \Z} \eps_y \eta_{2n-1}(y)} \big| \mathcal{F} \big] - \prod_{y \in \Z} \E \big[e^{it\eps_y \eta_{2n-1}(y)} \big| \mathcal{F} \big] \Big| e^{-t^2 d_{n,3}^2/2} dt
\]
and by Lemma \ref{lem2} and the  $\mathcal{F}$-measurability of the $\eta$'s:
\beqt
J_n &\leq& \frac{1}{2}  \sum_{x \neq y} \eta_{2n-1}(x) \eta_{2n-1}(y) \E [\eps_x \eps_y] \int_{|t| \leq b_n} t^2 dt.\\
&\leq& \frac{b_n^3}{6}
\sum_{x \neq y} \eta_{2n-1}(x) \eta_{2n-1}(y) \E [\eps_x \eps_y] .
\eeqt
Using the positivity of the correlations and the fact that we only work on $A_n$, we rewrite:
\beqt
J_n&\leq& \frac{n^{3 \delta_2}}{6n^{\frac{3}{2}+3\delta_3}} \; n^{\frac{1}{2}+\delta_2} \sum_{y=-2n}^{2n}  \eta_{2n-1}(y) \sum_{x=-2n,x \neq y}^{2n} \E [\eps_x \eps_y].\\
\eeqt
By stationarity of the associated r.v.'s $\eps_y$, we have for all $y \in [-2n,2n]$,
\[
\sum_{x=-2n,x \neq y}^{2n} \E [\eps_x \eps_y]= \sum_{x=-2n,x \neq y}^{2n} \E [\eps_0 \eps_{x-y}]=\sum_{z=-2n-y,z \neq 0}^{2n - y} \E [\eps_0 \eps_z] \leq \sum_{z=-4n, z \neq 0}^{4n} \E [\eps_0 \eps_z].
\]
Thus, still by stationarity,
\[
J_n \leq \frac{n^{3 \delta_2}}{6n^{\frac{3}{2}+3\delta_3}} \cdot n^{\frac{1}{2}+\delta_2} \cdot 2n \cdot 2 \sum_{z=1}^{4n} \E [\eps_0 \eps_z]
\]
In our case of summable power law decay of correlation, we have with $\alpha > 1$
\[
\E \big[\eps_0 \eps_y]=\mathcal{O} \big(|y|^{-\alpha})\; \Longrightarrow
\; \sum_{y=1}^{4n} \E [\eps_0 \eps_y]= \mathcal{O} (n^{1-\alpha})
\]
and thus
\[
J_n = \mathcal{O} \big(n^{-\beta}\big)
\]
with $\beta =3 \delta_3 + \alpha -1 - 4 \delta_2$.
\epr

Now, using (\ref{eqn4}), write with the usual notation $d_{n,3}= n^{\frac{1}{2}+\delta_3}$:
\[
\pee[A_n \setminus B_n | \mathcal{F}] \leq C  d_{n,3} (|I_n^0| +  |I_n^1 - I_n^0 | + |I_n^2|).
\]
Consider $\delta_3 > 2 \delta_2$. By the previous lemmata,   we have
\[
d_{n,3} \cdot |I_n^0|= \mathcal{O} \big(n^{-\frac{1}{4}+\delta_3 + \frac{\delta_1}{2}} \big), \; d_{n,3} \cdot |I_n^2|= \mathcal{O} \big(e^{-n^{2 \delta_2} / 2})
\]
and
\[
d_{n,3} \cdot |I_n^1 - I_n^0 |=\mathcal{O} \big(n^{\frac{1}{2}+\delta_3 - \beta} \big).
\]
To find a suitable $\delta>0$ such that Proposition \ref{prop1} holds, we need the following relations to be verified:
\begin{itemize}
\item $\delta_3 < \frac{1}{4}-\frac{\delta_1}{2}$.
\item $\frac{1}{2}+\delta_3 - \beta<0$,  or equivalently $\delta_3>2\delta_2 +\frac{1}{2}(\frac{3}{2} -\alpha)$
\end{itemize}
and we still need $\delta_3>2\delta_2$. We distinguish two cases:
\begin{itemize}
\item{$\alpha \in ]1,\frac{3}{2}[$}: the system reduces to
\[
\left\{
\begin{array}{lll} \delta_3>2\delta_2 +\frac{1}{2}(\frac{3}{2} -\alpha)\\
\\
\delta_3 < \frac{1}{4}-\frac{\delta_1}{2}.
\end{array}
\right.
\]
where $\delta_1$ and $\delta_2$ can be taken as small as possible so the existence of $\delta>0$ in Proposition \ref{prop1} requires that
\[
\frac{1}{2} \big( \frac{3}{2}-\alpha \big) < \frac{1}{4}
\]
i.e. $\alpha >1$, which is always verified under our hypothesis.
\item{$\alpha \geq \frac{3}{2}$}: the system reduces to
\[
\left\{
\begin{array}{lll} \delta_3>2\delta_2 \\
\\
\delta_3 < \frac{1}{4}-\frac{\delta_1}{2}.
\end{array}
\right.
\]
and one only has to choose $\delta_1$ and $\delta_2$ such that $2\delta_2 < \frac{1}{4} -\frac{\delta_1}{2}$ to find a suitable $\delta>0$.
\end{itemize}
This proves Proposition \ref{prop1}.

\epr

Combining Equations  (\ref{pn}), (\ref{srw}), (\ref{sqrtln}) and (\ref{eqn3}),  we obtain (\ref{eqn2}) and then (\ref{eqn1}). By Borel-Cantelli's Lemma, we get
\[
\pee[M_{T_n}=(0,0)\; \rm{i.o.}]=\pee[\pee[M_{T_n}=(0,0)\; \rm{i.o.}|\mathcal{G}]]=0
\]
and thus $(M_{T_n})_{n \in \N}$ is transient for $\nu$-a.e. orientation $\eps$. Theorem \ref{thm1} follows from Lemma \ref{lem1}.

\subsection{Proof of the functional limit theorem.}
\bp\label{pr1}
The sequence of random processes $n^{-3/4}(X_{[nt]})_{t\ge 0}$ weakly converges in the space ${\cal D}([0,\infty[,\R)$ to the process $(m\Delta_{t})_{t\ge 0}$.
\ep
\bpr
Let us first prove that the finite dimensional distributions of $n^{-3/4}(X_{[nt]})_{t\ge 0}$ converge to those of $(m\Delta_{t})_{t\ge 0}$ as $n\rightarrow\infty$.
We can rewrite for every $n\in\N$,
$$X_{n}=X_{n}^{(1)}+X_{n}^{(2)}$$
where
$$X_{n}^{(1)}=\sum_{y\in\Z} \epsilon_{y}\Big(\sum_{i=1}^{\eta_{n-1}(y)}\xi_{i}^{(y)} - m\Big)$$
and
$$X_{n}^{(2)}=m\sum_{y\in\Z} \eps_{y}\eta_{n-1}(y).$$
\bl\label{lem8}
The sequence of random variables $n^{-3/4} (X_{n}^{(1)})_{n\in\N}$ converges in
probability to 0 as $n\rightarrow +\infty$.
\el
\bpr
It is enough to prove the convergence to 0 for the $L^2$-norm.
$$
\E\Big[(X_{n}^{(1)})^{2}\Big]=\E\Big[\sum_{x,y\in\Z}\eps_{x}\eps_{y}
\sum_{i=1}^{\eta_{n-1}(x)}\sum_{j=1}^{\eta_{n-1}(y)}\E[(\xi_{i}^{(x)}-m)(\xi_{j}^{(y)}-m)|\mathcal{F} \vee \mathcal{G}]\Big]
$$
Since by independence of the $\xi_i^{(y)}$'s with both the vertical walk and the orientations, $$\E[(\xi_{i}^{(x)}-m)(\xi_{j}^{(y)}-m)|\mathcal{F} \vee \mathcal{G}]=\E[(\xi_{i}^{(x)}-m)(\xi_{j}^{(y)}-m)]m^2\delta_{i,j}\delta_{x,y},$$ we obtain
$$
n^{-3/2}\E\Big[(X_{n}^{(1)})^{2}\Big]=m^2 n^{-3/2}\sum_{x\in\Z} \eta_{n-1}(x)=m^2n^{-1/2}=o(1).
$$
\epr
\bl\label{lem9}
The finite dimensional distributions of $(n^{-3/4} X_{[nt]}^{(2)})_{t\ge 0}$ converge to those of $(m\Delta_{t})_{t\ge 0}$ as $n\rightarrow 0$.
\el
\bpr
Let $0 \leq t_{1} \leq t_{2} \leq \ldots \leq t_{k}$ and $\theta_{1},\theta_{2},\ldots,\theta_{k}\in \R$. By the definition of $X_{n}^{(2)}$, we have
$$n^{-3/4}\sum_{j=1}^{k}\theta_{j}X_{[nt_{j}]}^{(2)}=m n^{-\frac{3}{4}}
\sum_{j=1}^{k}\theta_{j}\sum_{y\in\Z}\eps_{y}\eta_{[nt_{j}]-1}(y).$$
For $\delta_{1}>0$, we define the event
$$D_{n}=\{\om\in\Omega;\ \max_{y\in\Z} \eta_{n}(y)< n^{\frac{1}{2}+\delta_1}\}.$$
One has
\beqt
\E\Big[\exp\Big(in^{-3/4}\sum_{j=1}^{k}\theta_{j}X_{[nt_{j}]}^{(2)}\Big)\Big]&=&\E\Big[
\E[\exp(imn^{-3/4}\sum_{j=1}^{k}\theta_{j}\sum_{y\in\Z}\eps_{y}\eta_{[nt_{j}]-1}(y))|{\cal F}]\Big]\\
&=&\E\Big[\E[\exp(imn^{-3/4}\sum_{j=1}^{k}\theta_{j}\sum_{y\in\Z}\eps_{y}\eta_{[nt_{j}]-1}(y))|{\cal F}]{\bf1}_{D_n}\Big]\\
&+&\E\Big[\E[\exp(imn^{-3/4}\sum_{j=1}^{k}\theta_{j}\sum_{y\in\Z}\eps_{y}\eta_{[nt_{j}]-1}(y))|{\cal F}]{\bf1}_{D_n^{c}}\Big]\\
&=&\Sigma_{1}(n)+\Sigma_{2}(n)\ \  (\mbox{say})
\eeqt
Firstly, we can use standard properties of the local time for the simple random walk on the line to estimate the second term:
$$|\Sigma_{2}(n)|\le \pee(D_n^{c})\le e^{-c n^{\delta_{2}}}$$
for some $c$ and $\delta_{2}$ strictly positive.\\*
Secondly, we compare on the particular set $D_n$ (on which uniformly in $y\in\Z$, the local time of the simple random walk is dominated by $n^{3/4}$) the characteristic function of the linear combinations of our process conditionally to the random walk with the marginal characteristic functions, using Lemma \ref{lem2}. Therefore we decompose
$$\Sigma_{1}(n)=\Sigma_{1,1}(n)+\Sigma_{1,2}(n)$$
where \beqt \Sigma_{1,1}(n)&=&\E\Big[{\bf1}_{D_n}\Big\{\E[\exp(im
n^{-3/4}\sum_{j=1}^{k}\theta_{j}
\sum_{y\in\Z}\eps_{y}\eta_{[nt_{j}]-1}(y))|{\cal F}]\\
&-&\prod_{y\in\Z}\E[\exp(im n^{-3/4}\eps_{y}\sum_{j=1}^{k}
\theta_{j}\eta_{[nt_{j}]-1}(y))|{\cal F}]\Big\}\Big]\\
\eeqt
and
$$\Sigma_{1,2}(n)=\E\Big[{\bf1}_{D_n}\prod_{y\in\Z}\E[\exp(im n^{-3/4} \eps_{y}\sum_{j=1}^{k}
\theta_{j}\eta_{[nt_{j}]-1}(y))|{\cal F}]\Big].$$
From Proposition 1 in \cite{KS},
we have that
\beqt
\lim_{n\rightarrow\infty}\Sigma_{1,2}(n)&=&\lim_{n\rightarrow\infty}
\E\Big[\exp\Big(-\frac{m^2}{2}n^{-\frac{3}{2}}\sum_{y\in\Z}
(\sum_{j=1}^{k}\theta_{j}\eta_{[nt_{j}]-1}(y))^2\Big)\Big]\\
&=&\E\Big[\exp\Big(-\frac{m^2}{2}\int_{-\infty}^{\infty}(\sum_{j=1}^{k}
\theta_{j}L_{t_{j}}(x))^2dx\Big)\Big] \mbox{ by Lemma 6 in \cite{KS}}\\
&=&\E\Big[\exp\Big(im\sum_{j=1}^{k}\theta_{j}\Delta_{t_{j}}\Big)\Big], \mbox{ see Lemma 5 in \cite{KS}}.
\eeqt

It remains to prove that $\Sigma_{1,1}(n)$ tends to 0 as $n$ goes to infinity.
By Lemma \ref{lem2}, we have that
\beqt
|\Sigma_{1,1}(n)|&\le &\frac{m^2}{2n^{3/2}}\sum_{x\neq y} \sum_{i,j=1}^{k} \theta_{i}\theta_{j}\E\Big[\eps_{x}\eps_{y}\eta_{[nt_{i}]-1}(x)
\eta_{[nt_{j}]-1}(y){\bf1}_{D_n}\Big]
\eeqt
Using the fact that we work on $D_n$, there exists $C>0$ such that
$$|\Sigma_{1,1}(n)|\le C \frac{n^{\frac{3}{2}+\delta_{1}}}{n^{\frac{3}{2}}}
\sum_{z=1}^{[t_k n]} \E[\eps_{0}\eps_{z}].$$
From the hypothesis on the power-law decay of correlations, there exists $\gamma >0$ such that
$$\sum_{z=1}^{n}\E[\eps_{0}\eps_{z}]={\cal O}(n^{-\gamma}).$$
So it is enough to choose $\delta_{1}<\gamma$ in order to have $\Sigma_{1,1}(n)=o(1)$ as $n$ goes to infinity.
\epr

From Lemma \ref{lem8} and Lemma \ref{lem9}, we deduce the convergence of the finite dimensional distributions of $n^{-3/4} (X_{[nt]})_{t\ge 0}$ to those of $(m\Delta_{t})_{t\ge 0}$.

In order to prove the weak convergence of
$(n^{-3/4}X_{[nt]})_{t\geq 0}$ to $(m\Delta_{t})_{t\geq 0}$ in
${\cal D}([0,\infty),\R)$, it remains to prove the tightness of
the family $(n^{-3/4}X_{[nt]})_{t\geq 0, n \geq 1}$ in  ${\cal
D}([0,\infty),\R)$. By Theorem 15.6 from Billingsley (\cite{bil}),
it is enough to prove that there exists $C>0$ such that for all
$t_{1}\le t \le t_{2}\in[0,T], T<\infty,$ for all $n\geq 1$,
\be\label{pro}
\E\Big[|X_{[nt_{2}]}-X_{[nt]}||X_{[nt]}-X_{[nt_{1}]}|\Big]\leq
C|t_{2}-t_{1}|^{\frac{3}{2}}. \ee Let us estimate \beqt
\E\Big[|X_{[nt_{2}]}-X_{[nt]}|^2\Big]&\le&  2m^2
\sum_{x,y\in\Z}\E\Big[(\eta_{[nt_{2}]-1}(x)-\eta_{[nt]-1}(x))(\eta_{[nt_{2}]-1}(y)-
\eta_{[nt]-1}(y))\Big]\E[\eps_{x}\eps_{y}]\\
&=& 2m^2 \sum_{z=-2n}^{2n} \E[\eps_{0}\eps_{z}] \sum_{x=-n}^{n} \E\Big[\Big(\sum_{k=[nt]}^{[nt_{2}]-1}{\bf 1}_{Y_{k}=x}\Big)\Big(\sum_{l=[nt]}^{[nt_{2}]-1}{\bf 1}_{Y_{l}=x+z}\Big)\Big]\\
&=& 2m^2 \sum_{z=-2n}^{2n} \E[\eps_{0}\eps_{z}] \sum_{k,l=[nt]}^{[nt_{2}]-1}
\pee[Y_{k}-Y_{l}=z]\\
&=& 2m^2 \Big\{2\sum_{z=-2n}^{2n} \E[\eps_{0}\eps_{z}]\sum_{k,l=[nt];k<l}^{[nt_{2}]-1}
\pee[Y_{l-k}=z]+[nt_{2}]-[nt]\Big\}.
\eeqt
Now, it is well-known that when $(Y_{k})_{k\ge 0}$ is a simple random walk on $\Z$, the probabilities of transition from 0 to $z$ satisfy uniformly in $z\in\Z$,
$$\pee[Y_{n}=z]={\cal O} \Big(\frac{1}{\sqrt{n}}\Big)$$
which implies that
\beqt
\sum_{k,l=[nt]; k<l}^{[nt_{2}]-1}\pee[Y_{l-k}=z]&=&{\cal O}\Big(([nt_{2}]-1-[nt])^{3/2}\Big)\\
&=&{\cal O}\Big(n^{3/2}(t_{2}-t_{1})^{3/2}\Big).
\eeqt
Using the hypothesis on the power-law decay of correlations and their positivity,
$$\sum_{z=-\infty}^{\infty}\E[\eps_{0}\eps_{z}]<\infty.$$
So we deduce that there exists $C>0$ such that
$$\E[|X_{[nt_{2}]}-X_{[nt]}|^2]\le C n^{3/2} |t_{2}-t_{1}|^{3/2}.$$
By Cauchy-Schwarz inequality, we obtain that there exists $C'>0$ such that
\beqt
n^{-3/2}\E[|X_{[nt_{2}]}-X_{[nt]}||X_{[nt]}-X_{[nt_{1}]}|]&\leq & n^{-3/2}\E[|X_{[nt_{2}]}-X_{[nt]}|^2]^{1/2}\E[|X_{[nt]}-X_{[nt_{1}]}|^2]^{1/2}\\
&\le & C' |t_{2}-t_{1}|^{3/2}
\eeqt
so the tightness is proved.
\epr

Let us recall that $M_{T_{n}}=(X_{n},Y_{n})$ for every $n\ge 1$.
The sequence of random processes $n^{-3/4}(Y_{[nt]})_{t\ge 0}$ weakly converges in ${\cal D}([0,\infty[,\R)$ to 0, thus the sequence of $\R^2-$valued random processes $n^{-3/4}(M_{T_{[nt]}})_{t\ge 0}$ weakly converges in ${\cal D}([0,\infty[,\R^2)$ to the process $(m\Delta_{t},0)_{t\ge 0}$.
Theorem 2.4 follows from this remark and the next lemma
\bl
The sequence of random variables $(\frac{T_{n}}{n})_{n\ge 1}$ converge in probability to $1+m$ as $n\rightarrow ~+\infty.$
\el
\bpr
Let us remark that
$$T_{n}=n+\sum_{y\in\Z}\sum_{i=1}^{\eta_{n-1}(y)}(\xi_{i}^{(y)}-m)
+m\sum_{y\in\Z}\eta_{n-1}(y).$$
Now, 
\beqt
\E\Big[\Big(\sum_{y\in\Z}\sum_{i=1}^{\eta_{n-1}(y)}(\xi_{i}^{(y)}-m)\Big)^2\Big]
&=&\sum_{x,y\in\Z} \E\Big[\sum_{i=1}^{\eta_{n-1}(x)}\sum_{j=1}^{\eta_{n-1}(y)}(\xi_{i}^{(x)}-m)(\xi_{j}^{(y)}-m)\Big]\\
&=& \sum_{x,y\in\Z} \E\Big[\sum_{i=1}^{\eta_{n-1}(x)}
\sum_{j=1}^{\eta_{n-1}(y)}\E[(\xi_{i}^{(x)}-m)(\xi_{j}^{(y)}-m)|{\cal F}\vee{\cal G}]\Big]\\
&=& m^2 \sum_{x,y\in\Z} \E\Big[\sum_{i=1}^{\eta_{n-1}(x)}
\sum_{j=1}^{\eta_{n-1}(y)} \delta_{i,j}\delta_{x,y}\Big]\\
&=& m^2\sum_{x\in\Z} \E[\eta_{n-1}(x)]\\
&=& m^2 n= \mathcal{O} \Big( n^2 \Big).
\eeqt
From this calculation and the fact that $\sum_{x\in\Z}\eta_{n-1}(x) =n$, we deduce the lemma.
\epr

\section{Examples}
Our framework includes this of \cite{CP} where  i.i.d. orientations are considered but
 it also includes orientations whose joint distribution is not a product measure. Natural examples
 of non-product measures are given by {\em Gibbs measures in statistical mechanics} (\cite{G}).
To destroy the independence of the random variables
 $\eps_y$, a family of measurable functions
$\Phi=(\Phi_A)_{A \in \s}$ indexed by the set $\s$ of finite
subsets of $\Z$, called {\em potential}, is introduced. For all $\eps, \Phi_A(\eps)$
represents the interaction between the random variables
$(\eps_y)_{y \in A}$ . A translation-invariant
measure $\nu^{\beta}$ on $\{-1,+1\}$ is a {\em Gibbs measure at
inverse temperature $\beta >0$} for the interaction $\Phi$ when it is
an equilibrium state regard to some variational properties in
terms of thermodynamic functions. An
equivalent definition characterizes Gibbs measures in terms of
continuity of their conditional probabilities w.r.t. the outside
of finite sets, or via the well-known {\em DLR equation}. In
some domains of temperature, there could be more than one Gibbs
measure for an interaction $\Phi$, and we then say that a {\em
phase transition} holds. In this Gibbsian context, FKG property requires {\em ferromagnetic pair interactions} (\cite{FKG}), i.e. a potential $\Phi$ such that $\Phi_A = 0$ if
$\rm{card}(A) > 2$ and $\Phi_{\{i,j\}}(\eps)= J(i,j) \eps_i \eps_j$ with $J(i,j) \geq 0$ for all $i,j \in \Z$. This provides us a wide
family of examples suitable to our set-up.

\begin{enumerate}
\item{Ferromagnetic nearest neighbors Ising model}:\\
The coupling is  translation-invariant, positive for nearest neighbors pairs $\{i,j\} \subset \Z$ and null otherwise:
\begin{displaymath}
 J(i,j)\; = \;\left\{
\begin{array}{lll} J \; \geq \; 0  \; &\rm{if} \;  |i-j|=1&\\
\\
0 \; \; &\rm{otherwise.}&
\end{array}
\right.
\end{displaymath}
It is well known (\cite{G}) that this one-dimensional model does not exhibit a phase transition.
 By translation invariance of the potential, the unique Gibbs measure is then translation
  invariant, and equivalently the family $\eps$ is stationary. Moreover, the
  decay of correlations of the random field  $\eps$ of law $\nu$ is known
   to be exponential (\cite{cas,DS}). The absence of phase transition  and the
   translation invariance property of the potential prove also that the orientations
   are symmetric and thus fit in our framework.

\item{Long range ferromagnetic Ising model}:\\
It is similarly defined but the coupling $J(i,j)$ is non
 null for any pair $\{i,j\}$ and has a power law decay: there exists  $\alpha >1$ and $J \geq 0$ such that
\[
J(i,j)=J.|i-j|^{-\alpha}.
\]
Depending on the value of $\alpha$, there could be a phase transition in some domains of temperature, and in particular two different regimes with summable power law decay of correlations are relevant in our set-up (see \cite{AZ,D,Na,Rue}).
\begin{enumerate}
\item{$\alpha > 2$}.\\
There is no phase transition and the Gibbs measure
 for this potential  is translation invariant. The variables are thus
  symmetric and one could also learn in the literature
  that $\alpha$ also governs a power law decay of  correlations:
\[
\E \big[\eps_0 \eps_y]= \mathcal{O} \big( {\mid y \mid}^{-\alpha} \big) \; \rm{when} \;
\mid y \mid  \longrightarrow \; + \infty.
\]
\item{$\alpha \in \;]1,2[$}.\\
There exists a critical inverse temperature $\beta_c > 0$ which
separates the domain of temperature in two different regimes.
 In the
high temperature regime, there is no phase transition and the picture is as above: the power law decay is the same as this of the interaction, i.e. governed by $\alpha$. By translation-invariance of the unique phase, this provides us examples with very slow but summable power law decay of correlations for which our theorems hold.

The later is not necessarily true at low (or critical) temperature. Firstly, when a phase transition holds, the random variables are not necessarily centered and secondly, the power law decay of correlations  is governed by $\alpha-1 < 1$, and thus the correlations are not summable.
\end{enumerate}

\end{enumerate}

\section{Comments}
We have extended the results of \cite{CP} to positively
 correlated  orientations and solved one
  of their open problems. In particular, we have proved that the simple
   random walk is still transient for ferromagnetic models in absence of phase
   transition. As the walk can be  recurrent for deterministic orientations,
    it would be interesting to perturb deterministic cases in order to get a
    full picture of the transience versus recurrence properties and identify
    a sort of phase transition. As a perturbation of the alternate lattice
     on which the walk is recurrent, we are studying  anti-ferromagnetic
     systems (i.e. the case of a coupling $J \leq 0$) for which the recurrence
     might  be conserved at low temperature. Similarly, one could consider
     negatively correlated orientations but  this requires finer results on such distributions
      and a complete theory of negative dependence has not been established yet
(see e.g. \cite{pem,cam}).

The ingredients used to prove the functional limit theorems still hold for any $m > 0$ and Theorem \ref{thm2} should
therefore also be satisfied  for more  general random walks  than the
simple random walk on the lattice $\mathbb{L}^\eps$. This question is currently under
considerations.\\

{\bf Acknowledgments:}
We thank D. P\'etritis for having introduced us to this subject
and indicated to us the proof of Lemma \ref{lem6}. We are also
grateful to Aernout van Enter and Frank Redig for many useful comments and
references. This work has been partially achieved under the RDSES
program of the ESF.


\addcontentsline{toc}{section}{\bf References}
\begin{thebibliography}{14}

\bibitem{AZ} M. Aizenman and R. Fern\'andez. Critical exponents
for long range interactions. {\em Lett. Math. Phys.} {\bf 16}, no
1:39--49, 1988.

\bibitem{bil} P. Billingsley. Convergence of probability measures. Wiley, New-York, 1968.

\bibitem{cam} C. Cammarota. Positive and negative correlations for conditional Ising distributions.
{\em Rev. Math. Phys.} {\bf 14}, no 10:1099--1113, 2002.

\bibitem{CP} M. Campanino and D. P\'etritis. Random walks on randomly oriented lattices.
{\em Mark. Proc. Relat. Fields}, in press, 2003.
 Available at http://name.math.univ-rennes1.fr/dimitri.petritis/articles.html.

 \bibitem{CP2} M. Campanino and D. P\'etritis. On the physical relevance of random walks: an example of random walks on randomly oriented lattices, in "Random walks and geometry", in press, 2003.

\bibitem{cas} M. Cassandro and E. Olivieri. Renormalization group
and analycity in one dimension: a proof of Dobrushin's theorem. {\em
Commun. Math. Phys.} {\bf 80}, no 2:255--269, 1981.

\bibitem{DS} R. L. Dobrushin and S. B. Shlosman. Completely analytical interactions:
 constructive description. {\em J. Statist. Phys.} {\bf 46}, no 6:983--1014, 1987.

\bibitem{D} F. J. Dyson. Existence of a phase-transition in a one-dimensional Ising ferromagnet.
 {\em Commun. Math. Phys.} {\bf 12}:91--107, 1969.
\bibitem{FFS} R. Fern\'andez, J. Fr\"ohlich and  A. D. Sokal.
 {\em Random-Walks, Critical Phenomena, and Triviality in Quantum Field Theory}.
Springer-Verlag (Texts and Monographs in Physics), 1992.

\bibitem{FKG} C. M. Fortuin, J. Ginibre and P. W. Kasteleyn.
Correlation inequalities on some partially ordered sets.
 {\em Commun. Math. Phys.} {\bf 22}:89--103, 1971.


\bibitem{G} H. O.~Georgii.  \newblock Gibbs Measures and Phase
  Transitions. \newblock Walter de Gruyter (de Gruyter Studies in
  Mathematics, Vol.\ 9), Berlin-New York, 1988.





\bibitem{KS} H. Kesten and F. Spitzer. A limit theorem related to a new class of self similar
processes. {\em Z. Wahrsch. Verw. Gebiete} {\bf 50}:5--25, 1979.

\bibitem{LR} P. Leroux. Coassociate grammars, periodic orbits and quantum random walks over $\Z^d$. Pr\'epublication Irmar,
Universit\'e de Rennes 1, 2003.

\bibitem{Na} A. Naimzhanov. Weakening of correlations in one-dimensional systems with long-range interaction. {\em Theor. Math. Phys.} {\bf 38}:64--68, 1979.

\bibitem{NW} C. M. Newman and A. L. Wright. An invariance principle for certain dependent
sequences. {\em Ann. Probab.} {\bf 9}, no 4:671--675, 1981.

\bibitem{pem} R. Pemantle. Towards a theory of negative dependence. Probabilistic techniques in equilibrium and nonequilibrium statistical physics. {\em J. Math. Phys.} {\bf 41}, no. 3:1371--1390, 2000.

\bibitem{Rue} D. Ruelle. Statistical mechanics of a one-dimensional lattice gas. {\em Commun. Math. Phys.} {\bf 9}:267--278, 1968.
\end{thebibliography}
\end{document}